\documentclass[a4paper,12pt] {elsart}
\usepackage[body={165mm,225mm}]{geometry}
\usepackage[latin1]{inputenc}
\usepackage{amsfonts,amssymb}

\usepackage{graphicx}

\usepackage{amsmath}

\newtheorem{proposition}{Proposition}

\def\CaixaPreta{\vrule Depth0pt height5pt width5pt}
\def\comecaprova{\noindent {\bf Proof:} \hspace{3mm}}
\def\terminaprova{\hfill \CaixaPreta \vspace{5mm}}
\def\Z{\mathbb Z}
\def\R{\mathbb R}
\def\S{\mathbb S}
\def\P{\mathbb P}

\begin{document}
\begin{frontmatter}

\title{An Affine Linear Solution for the 2-Face Colorable Gauss Code Problem
in the Klein Bottle and a Quadratic System for Arbitrary Closed
Surfaces}
\author[A1]{S\'ostenes Lins}
\author[A2]{\hspace{-1mm}, Emerson Oliveira-Lima}
\author[A3]{\hspace{-1mm}, Valdenberg Silva}
\address[A1]{Dept. Matem\'atica da UFPE - Recife - Brazil}
\address[A2]{Dept. Matem\'atica da UNICAP - Recife - Brazil}
\address[A3]{Dept. Matem\'atica da UFSE - Aracaju - Brazil}

\bibliographystyle{alpha}
\maketitle

\begin{abstract}

Let $\overline{P}$ be a sequence of length $2n$ in which each
element of $\{1,2,\ldots,n\}$ occurs twice. Let $P'$ be a closed
curve in a closed surface $S$ having $n$ points of simple
auto-intersections, inducing a 4-regular graph embedded in $S$
which is 2-face colorable. If the sequence of auto-intersections
along $P'$ is given by $\overline{P}$, we say that is a {\em $P'$
$2$-face colorable solution for the Gauss Code $\overline{P}$ on
surface $S$} or a {\em lacet for $\overline{P}$ on $S$}. In this
paper we present a necessary and sufficient condition yielding
these solutions when $S$ is Klein bottle. The condition take the
form of a system of $m$ linear equations in $2n$ variables over
$\Z_2$, where $m \le n(n-1)/2$. Our solution generalize solutions
for the projective plane and on the sphere. In a strong way, the
Klein bottle is an extremal case admitting an affine linear
solution: we show that the similar problem on the torus and on
surfaces of higher connectivity are modelled by a quadratic system
of equations.

\begin{description}\item [Keywords:] {Gauss code problem, lacets, closed surfaces, $4$-regular graphs,
medial maps (of graphs on surfaces), face colorability}

\end{description}

\begin{description}\item [Mathematical Subject Classification:] {primary 05C10; secondary 05C85}

\end{description}

\end{abstract}
\end{frontmatter}

\section{Introduction}
A Gauss code $\overline{P}$ is a cyclic sequence in the set of
labels $E = \{1,2,\ldots,n\}$ in which each $x \in E$ occurs
twice. Let $P'$ be a closed curve in a closed surface $S$ having
$n$ points of simple auto-intersections, inducing a 4-regular
graph embeddded into $S$ such that the cyclic sequence of
auto-intersections reproduces $\overline{P}$. We say that $P'$
realizes $\overline{P}$ in $S$ and it is called a {\em lacet} for
$\overline{P}$. If the embedding of $P'$ produces a $2$-face
colorable map, then $P'$ is called a {\em $2$-colorable lacet} for
$\overline{P}$. Without the $2$-face colorability condition the
algebra derived from maps with a single zigzag \cite{Lins1980} is
not available and an entirely different problem arises. Here we
consider only $2$-colorable lacets. In this case, $P'$ is the {\em
medial map} (\cite{GodsilRoyle2001}) of a map $M$ formed by a
graph $G_M$ embedded into $S$. $M$ has a single zigzag
\cite{Lins1982}. The dual of $M$ is denoted $D$ and its phial
(\cite{Lins1982}) is denoted $P$. $G_P$, the graph of $P$, has a
single vertex (corresponding to the single zigzag). Previous work
on the Gauss code problem can be found in \cite{Shank1975},
\cite{Rosenstiehl1976a}, \cite{Rosenstiehl1976b},
\cite{LovazMax1976}, \cite{Lins1980} and
\cite{LinsRichterShank1987}. The last two works solve the $2$-face
colorable problem for the case of the projective plane. The
previous works deal with the planar case in which the $2$-face
colorability is granted. In the present work we algorithmically
solve the problem for the Klein bottle.

This problem has been recently tackled by in
\cite{CrapoRosenstiehl2001}.  In this paper they introduce the
terminology {\em lacet} meaning $2$-colorable lacets and develop a
theory which has its origin in the basic algebraic fact appearing
in \cite{Lins1980} and in \cite{LinsRichterShank1987} connecting
the surface $S$ and the intersection of the cycle spaces of $G_M$
and $G_D$. They find conditions for their realization in the torus
(Theorem 19) and in the Klein bottle (Theorem 20) in terms of the
the existence of a pair of $0-1$ vectors with certain properties.
However, the verification of the existence of these vectors is
left undiscussed. In fact, to verify their existence leads to an
exponential number of trials. So the theorems do not provide
polynomial algorithms of even good characterizations (in the sense
of Edmonds) for realizability. In this technical sense the
problems are not solved in \cite{CrapoRosenstiehl2001}.

For the the case of the Klein bottle it is solved here. It remains
open for the case of $\S^1 \times \S^1$. In a strong way, the
Klein bottle is an extremal case admitting an affine linear
solution: we show that the similar problem on the torus and on
surfaces of higher connectivity are modelled by a quadratic system
of equations. We also provide a way for finding the smallest
connectivity of a surface realizing a Gauss code as a lacet in
terms of deciding whether or not a quadratic system of equations
over $\Z_2$ is consistent. We also provide a way for finding the
smallest connectivity of a surface realizing a Gauss code as a
lacet in terms of deciding whether or not a quadratic system of
equations over $\Z_2$ is consistent.

The paper is organized as follows. In Section 2 we show an example
to help the reader in understanding our definitions and
motivation. This example is used throughout the paper. In Section
3 we give the statement of the main results. Section 4 briefly
reviews the theory of Combinatorial Maps as given in
\cite{Lins1982} to prove the {\em Parity Theorem} \cite{Lins1980},
which is needed in this work as a Lemma. In Section 5 we prove
some basic lemmas and the main results. Section 6 deals with
2-face colorable lacets in general surfaces. We prove that this
general problem is modelled as a quadratic system of equations
over $\Z_2$. Finally, short Sections 7 and 8 are concluding
remarks and acknowledgements.

\section{An example and the algebraic tools}
Consider the example of a kleinian map $M$ with a single zigzag
given on the left of Fig. 1. The cyclic sequence of edges visited
in the zigzag is $$P=(1,4,5,6,5,4,3,8,7,3,2,-1,-2,8,7,-6).$$ This
can be easily followed in the $2$-colorable lacet $P'$. We think
of $P'$ as the medial map of $M$. The direction of the first
occurrence of an edge of $G_M$ defines the orientation. Edges
$3$,$4$,$5$,$7$,$8$ are traversed twice in the positive direction
(they correspond to black circles in the medial map) and edges
$1$,$2$,$6$ are traversed once in the positive direction and once
in the negative direction (they correspond to white circles). The
reason for the notation $P$ is that the signed cyclic sequence
defines the phial map $P$ (whence also $M$, $D$ and $P^\sim,$ as
well as the surface of $M$) and vice-versa, the phial defines the
sequence. For the algebraic concepts we refer to
\cite{Godement1968}. For the graph terminology to
\cite{BondyMurty1976} and to \cite{GodsilRoyle2001}. For more
background on graphs embedded into surfaces we refer to
\cite{Giblin1977}. Given a map $M$ with a single $z$-gon we define
linear functions (over the field $\Z_2$) $i_{\overline{P}}: 2^E
\rightarrow 2^E$ and $\kappa_P: 2^E \rightarrow 2^E$ as follows.
They are defined in the singletons and extended by linearity. Let
$i_{\overline{P}}(\{x\})$ be the set of edges occurring once in
the cyclic sequence $P$ between the two occurrences of $x$. Let
$\kappa_P(x)=\{x\}$ if $x$ is traversed twice in the same
direction in the zigzag path ($x$ is a black vertex in the medial
map), and $\kappa_P(x)=\emptyset$, if $x$ is traversed in opposite
direction in the zigzag path ($x$ is a white vertex in the medial
map). Let $c_P = \kappa_P + i_{\overline{P}}$ and ${\mathcal
V^\perp}$, ${\mathcal F^\perp}$ the cycle spaces of $G_M$ and
$G_D$, respectively. It is easy to verify that $c_P(x)$ is the set
of edges occurring once in a closed path in $G_M$. Therefore,
$c_P(x) \in {\mathcal V}^\perp$. In the above figure we see that
$c_P(1)=\emptyset \cup \{2,8,7,6\}, c_P(3)=\{3\} \cup \{8,7\}$
and, indeed, $\{2,8,7,6\}$ and $\{3,8,7\}$ are in ${\mathcal
V}^\perp$. From the definitions, it follows that if $P$ has a
single vertex, for any $x$, $\kappa_P(x) + \kappa_{P^\sim}(x) =
\{x\}$ and that $c_{P^\sim}(x)+c_P(x)=\{x\}$ and so,
$c_{P^\sim}+c_P$ is the identity linear transformation.

%-----------------------------------
    \vspace{0.5cm}
    \begin{figure}%[!h]
        \begin{center}
        \includegraphics{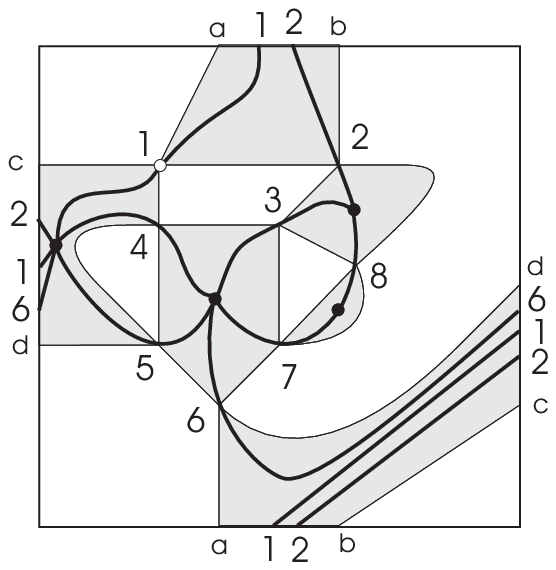} \hspace{20mm}\includegraphics{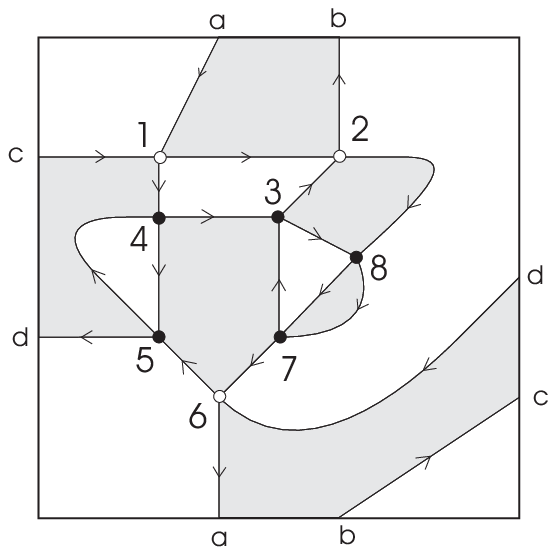}\\
            \caption{{\sf A kleinian map $M$ with $|V(P)|=1$ and its medial
            $P'$. Gauss code $\overline{P}=1456543873212876$}}
            \label{fig:Kuratowski14}
        \end{center}
    \end{figure}
%-----------------------------------

Let $b_P = c_{P^\sim} \circ c_P = c_P \circ c_{P^\sim}.$ The image
of $b_P$ is  ${\mathcal V}^\perp \cap {\mathcal F}^\perp$.
Moreover the dimension of this subspace is the connectivity of
$S$. These facts were first proved in \cite{Lins1980}. They  also
appear in \cite{LinsRichterShank1987}, in
\cite{CrapoRosenstiehl2001}. For the above example we have the
following values for functions $c_P$, $c_{P^\sim}$ and $b_P$:

\begin{center}
$\begin{array}{|l|l|l|} \hline c_P(1)=\{2,8,7,6\} &
c_{P^\sim}(1)=\{1,2,8,7,6\} & b_P(1)=\{2,4,5,6,7,8\}\\ \hline
c_P(2)=\{1\} & c_{P^\sim}(2)=\{2,1\} & b_P(2)=\{1,2,6,7,8\}\\
\hline c_P(3)=\{3,8,7\} & c_{P^\sim}(3)=\{8,7\} &
b_P(3)=\emptyset\\ \hline c_P(4)=\{4,6\} & c_{P^\sim}(4)=\{6\} &
b_P(4)=\{1,4,5\}\\ \hline c_P(5)=\{5,6\} & c_{P^\sim}(5)=\{6\} &
b_P(5)=\{1,4,5\}\\ \hline c_P(6)=\{1,4,5\} &
c_{P^\sim}(6)=\{6,1,4,5\} & b_P(6)=\{1,2,6,7,8\}\\ \hline
c_P(7)=\{7,3,1,8\} & c_{P^\sim}(7)=\{3,1,8\}&
b_P(7)=\{1,2,6,7,8\} \\ \hline c_P(8)=\{8,7,3,1\} &
c_{P^\sim}(8)=\{7,3,1\}& b_P(8)=\{1,2,6,7,8\} \\ \hline
\end{array}$
\end{center}

\section{Statement of the Main Results}

If $S$ is the Klein bottle, then dim(Im$(b_P))=2=\xi(S)$. Since
the underlying field is $\Z_2$, there are at most $4$ distinct
values in the image of $b_P$. Indeed we can be more specific. Let
${\mathcal O}=\{x \in E \ | \ |i_{\overline{P}}(x)|\ {\rm is \
odd}\}$ and ${\mathcal E}$ be the complementary subset of edges.

\begin{proposition} [Main Lemma] If $P'$ is the $2$-colored lacet for ${\overline{P}}$
which is in the $2$-sphere, $\S^2$, in the real projective plane,
$\R\P^2$, or in the Klein bottle, $\S^1
\times^{\raisebox{0.5mm}{\hspace{-3.2mm}$\sim$}} \S^1$, then there
are partitions $({\mathcal O}_0,{\mathcal O}_1)$ of ${\mathcal O}$
and $({\mathcal E}_0,{\mathcal E}_1)$ of ${\mathcal E}$ satisfying

\begin{center}
\begin{tabular}{lcrl}
$b_P(x)$&$=$&${\mathcal O}_0 \cup {\mathcal E}_1$,&if $x \in {\mathcal O}_0$\\
&$=$&${\mathcal O}_1 \cup {\mathcal E}_1$,&if $x \in {\mathcal O}_1$\\
&$=$&$\emptyset$,&if $x \in {\mathcal E}_0$\\
&$=$&${\mathcal O}$,&if $x \in {\mathcal E}_1$.\\
\end{tabular}
\end{center}

Moreover, if $P'$ is in the Klein bottle, then ${\mathcal O}_0
\neq \emptyset$ and at most one of ${\mathcal O}_1$ and ${\mathcal
E}_1$ is empty. \label{prop:teorema1}
\end{proposition}

This proposition is already proved for the cases of $\S^2$ and
$\R\P^2$: Lins$[1980]$, Lins, Richter and Shank$[1987]$. For the
Klein bottle the proof is given in Section 5.

For the example of Fig. 1 the parts are ${\mathcal E}_0 = \{3\}$,
${\mathcal E}_1 = \{1\}$, ${\mathcal O}_0 = \{2,6,7,8\}$,
${\mathcal O}_1 = \{4,5\}$.

\begin{proposition}[Main Theorem] Let $\overline{P}$ be a given Gauss
code with $E=\{1,2,\ldots,n\}$. For some $m \le n(n-1)/2$, there
exists an $(m \times 2n)$-matrix $L=L(\overline{P}) \in \Z^{m
\times \ 2n}_2$ (computable from $\overline{P}$) and a column
$2n$-vector $r=r(\overline{P})\in \Z^{2n}_2$ (also computable from
$\overline{P}$) such that the following characterization holds:
each $2$-face colorable realization of $\overline{P}$ in $\S^2$,
in $\R\P^2$ or in $\S^1
\times^{\raisebox{0.5mm}{\hspace{-3.2mm}$\sim$}} \S^1$ corresponds
in an $1-1$ way to a solution $\xi \in \Z^{2n}_2$ of the linear
system $L \xi = r$. If the code is not realizable in these
surfaces, then there exists an $m$-column vector $\nu \in
\Z^{m}_2$ such that $\nu^T L = 0$ and $\nu^T r = 1$.
\label{prop:MainTheorem}
\end{proposition}

The proof of the Main Theorem is also postponed to Section 5.

\section{Combinatorial maps and a Parity Theorem}

A {\em topological map} $M^t=(G,S)$ is an embedding of a graph $G$
into a closed surface $S$ such that $S\backslash G$ is a
collection of disjoint open disks, called {\em faces}. By going
around the boundary of a face and recalling the edges traversed we
define a {\em facial path} of $M^t$, which is a closed path in
$G$. Note that a facial path is obtained starting in an edge and
by choosing at each vertex always the rightmost or always the
leftmost possibility for the next edge. If we alternate the
choice, then the result is a {\em zigzag path}, or simply a {\em
zigzag}. Even if the surface is non-orientable these left-right
choices are well defined, because they are local. For more
background on graphs embedded into surfaces see \cite{Giblin1977}.
To make our objects less dependent of topology we use a
combinatorial counterpart for topological maps introduced in
\cite{Lins1982}. A {\em combinatorial map} or simply a {\em map}
$M$ is an ordered triple $(C_M,v_M,f_M)$ where: $(i)$ $C_M$ is a
connected finite cubic graph; $(ii)$ $v_M$ and $f_M$ are disjoint
perfect matchings in $C_M$, such that each component of the
subgraph of $C_M$ induced by $v_M\cup f_M$ is a {\em polygon}
(i.e. a non-empty connected subgraph with all the vertices having
two incident edges) with $4$ edges and it is called an {\em
$M$-square}.

From the above definition, it follows that $C_M$ may contain
double edges but not loops. A third perfect matching in $C_M$ is
$E(C_M)-(v_M\cup f_M)$ and is denoted by $a_M$. The set of
diagonals of the $M$-squares, denoted by $z_M$, is a perfect
matching in the complement of $C_M$. The edges in
$v_M,f_M,z_M,a_M$ are called respectively $v_M$-edges,
$f_M$-edges, $z_M$-edges, $a_M$-edges. The graph $C_M \cup z_M$ is
denoted by $Q_M$, and is a regular graph of valence 4. A component
induced by $a_M \cup v_M$ is a polygon with an even number of
vertices and it is called a $v$-gon. Similarly, we define an
$f$-gon, and a $z$-gon, by replacing $v$ for $f$ and $v$ for $z$.
Clearly, the $f$-gons and $z$-gons of $C_M$ correspond to the
facial paths and the zigzags of $M^t$. To avoid the use of colors
the $M$-squares are presented in the pictures as rectangles in
which the short sides $(s)$ are $v_M$-edges, the long sides
$(\ell)$ are $f_M$-edges and the diagonals $(d)$ are $z_M$ edges.
An {\em $M$-rectangle with diagonals} or simply an {\em
$M$-rectangle} (being understood that the diagonals are present)
is a component induced by $v_M, f_M, z_M$. The set of
$M$-rectangles is denoted by $R$. If $\pi$ is a permutation of the
symbols $s\ell d$, and $R'\subseteq R$ subset of rectangles of
$M$, then $M(R':\pi)$ denotes the map obtained from $M$ by
permuting the short sides, the long sides and the diagonals
according to $\pi$ in all $r \in R'$. Let $M(r: \pi)$ denote
$M(\{r\}: \pi)$. The {\em dual map} of $M$ is the map $D=M(R:\ell
sd)$; $D$ and $M$ have the same $z$-gons and the $v$-gons and
$f$-gons interchanged. The {\em phial map} of $M$ is the map
$P=M(R:d\ell s)$; $P$ and $M$ have the same $f$-gons and the
$v$-gons and $z$-gons interchanged. The {\em antimap} of $M$ is
the map $M^\sim=M(R:sd\ell)$; $M$ and $M^\sim$ have the same
$v$-gons and the $f$-gons and $z$-gons interchanged. The pairs
$(M,D)$, $(M,P)$, $(M,M^\sim)$ constitute the {\em map dualities}
introduced in $[{\rm Lins}, 1982]$. The dual of $P$ is $D^\sim$
and the dual of $M^\sim$ is $P^\sim$. Let
$\Omega(M)=\{M,D,P,M^\sim,D^\sim,P^\sim\}$ and
$\Omega^\star(M)=\{M(R':\pi)  \ | \ R' \subseteq R, \pi {\rm \
permutation\ of\ } s\ell d \}$. Note that $\Omega(M) \subseteq
\Omega^\star(M)$ and that any member of $\Omega^\star(M)$ has $R$
as its set of rectangles.

%-----------------------------------
    \vspace{0.5cm}
    \begin{figure}%[!h]
        \begin{center}
        \includegraphics{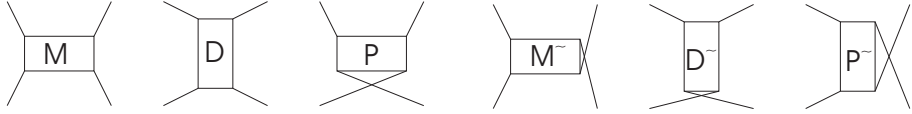} \\
            \caption{\sf How a neighborhood of each rectangle is modified in
                      the members of $\Omega(M)$}
            \label{fig:OmegaM}
        \end{center}
    \end{figure}
%-----------------------------------

Given a map $M$ and its dual $D$, there exists a closed surface,
denoted by ${\rm Surf}(M,D)$ where $C_M=C_D$ naturally embeds.
Consider the $v$-gons, the $f$-gons and the $M$-squares bounding
disjoint closed disks. Each edge of $C_M$ occurs twice in the
boundary of this collection of disks. Identify the collection of
disks along the two occurrences of each edge. The result is a
closed surface and $C_M$ is {\em faithfully embedded on it},
meaning that the boundaries of the faces are {\em bi}colored
poly\hspace{-0.8mm}{\em{gons}} or {\em bigons}. Similarly, there
are surfaces ${\rm Surf}(D^\sim,P)$ and ${\rm
Surf}(P^\sim,M^\sim)$.

We define a function $\psi$ which turns out to be a bijection from
the set of maps onto the set of $t$-maps. We denote $\psi(M)$ by
$M^t$. Given a map $M$, to obtain $M^t$ we proceed as follows.
Consider the $t$-map $(C_M,S)$, where $S={\rm Surf}(M,D)$, given
by the faithful embedding of $M$. The $v$-gons, the $f$-gons and
the $M$-squares are boundaries of (closed, in this case) disks
embedded (and forming) the surface $S(M)$. Shrink to a point the
disjoint closed disks bounded by $v$-gons. The $M$-squares, then,
become bounding digons. Shrink each such bounding digon to a line,
maintaining unaffected its vertices. With these contractions,
effected in $S$, $t$-map $(C_M,S)$ becomes, by definition,
$M^t=(G_M,S)$. Graph $G_M$ is called the {\it graph induced} by
$M$. A combinatorial description of $G_M$ can be given as follows:
the vertices of $G_M$ are the $v$-gons of $M$; its edges are the
squares of $M$; the two ends of an edge of $G_M$ are the two
$v$-gons (which may coincide and the edge is a {\em loop}) that
contain the $v_M$-edges of the corresponding $M$-square. It is
evident that $\psi$ is inversible: given a $t$-map we replace each
edge by a bounding digon in its surface, and then expand each
vertex to a disc in order to obtain a cellular embedding of a
cubic graph. Therefore, $\psi^{-1}$ is well-defined; in fact, it
is the dual of a useful construction in topology, namely,
barycentric division. Thus, $\psi$ is a bijection from the set of
maps onto the set of $t$-maps. It can be observed that $\psi$
induces a bijection from the set of $M$-rectangles onto the set of
edges of $G_M$. We use this bijection to identify the sets $R$ and
$E(G_M)$. Via $R$, which is invariant for the members of
$\Omega^\star(M)$, we identify $E(G_M)$ and $E(G_{M'})$ for $M'
\in \Omega^\star(M)$. Denote these identified sets of edges by
$E$.

Consider the function $\psi^M_c$, from the cycle space of
$C_M,CS(C_M)$, onto the cycle space of $G_M,CS(G_M)={\mathcal
V}^\perp$. It is defined as follows: for $S\in CS(C_M)$, an edge
$s\in E$ is in $\psi^M_c(S)$ if each square $s\in SQ(M)$ meets $S$
in exactly one $f_M$-edge. With this definition, it is evident
that $\psi^M_c(S)$ is a cycle in $G_M$ and that $\psi^M_c$ is
surjective.

\begin{proposition}[Lins 1980]
 $\psi^M_c$ is a homomorphism. Its kernel
is the subspace of $CS(C_M)$ generated by the $v$-gons and the
squares of $M$. \label{prop:lema1} \end{proposition}

\comecaprova Let $S_1$ and $S_2$ be cycles in $C_M$. We must show
that $\psi^M_c(S_1+S_2)=\psi^M_c(S_1)+\psi^M_c(S_2)$. An edge $e$
of $G_M$ is in $\psi^M_c(S_1+S_2)$ if, and only if, exactly one
$f_M$-edge of square $e\in SQ(M)$ is in $S_1+S_2$. Therefore,
$e\in\psi^M_c(S_1+S_2)$ if, and only if, an even number of
$f_M$-edges of square $e$ belongs to a fixed member of the set
$\{S_1,S_2\}$, and one $f_M$-edge of square $e$ belongs to the
other. The latter statement is equivalent to
$e\in\psi^M_c(S_1)+\psi^M_c(S_2)$. This proves that $\psi^M_c$ is
homomorphism. The image under $\psi^M_c$ of any square or any
$v$-gon of $M$ is the null cycle in $G_M$. Thus, the space
generated by the squares and $v$-gons is contained in ${\rm
Ker}(\psi^M_c)$. Conversely, suppose that $S\in{\rm
Ker}(\psi^M_c)$. The intersection of $S$ with the edges of an
arbitrary square has zero or two $f_M$-edges. Denote by $T$ the
cycle formed by the union of the squares $Q\in SQ(M)$ such that
$eQ\cap S$ contains the two $f_M$-edges of $Q$. It follows that
$S+T$ has no $f_M$-edges. Since $S+T$ is the edge set of a
collection of polygons, ($C_M$ is cubic), it follows that $S+T$ is
formed by a collection of $v$-gons whose edge set we denote by
$U$. Hence, $S=T+U$, with $T$ induced by squares, $U$ induced by
$v$-gons. Therefore, we conclude that ${\rm Ker}(\psi^M_c)$ is
contained in the space generated by the edge sets of squares and
$v$-gons of $M$. The proof is complete. \terminaprova

Since an element of the kernel of $\psi^M_c$ has an even number of
edges of $C_M$, it follows that if $\psi^M_c(S_1)=\psi^M_c(S_2)$,
then $|S_1|\equiv|S_2|$ $\mod\,2$. This observation makes the
following definition meaningful. A cycle $S$ in $G_M$ is called an
{\it $r$-cycle} in $M^t$ if $\psi^M_c(S')=S$ and $|S'|$ is odd,
for some cycle $S'$ in $C_M$. If $|S'|$ is even and
$\psi^M_c(S')=S$, then we say that $S$ is an {\it $s$-cycle} in
$M^t$. {\em $r$-Circuits} are minimal $r$-cycles. We observe that
the $r$-circuits in $M^t$ are precisely the orientation-reversing
polygons in $M^t$. This topological notion is not used; we work
with our parity definition of $r$-cycle. Observe that a subset
$T\subseteq E$ is a boundary in $M^t$ if, and only if, there
exists a cycle $T'$ of $C_M$, such that $\psi^M_c(T')=T$, and $T'$
can be written as the mod $2$ sum of nome subsets of $v$-gons,
$f$-gons, and squares. For a map $M$, two cycles in $G_M$ are {\em
homologous mod $2$} if their symmetric difference is a boundary in
$M^t$. Homology mod $2$ is, thus, an equivalence relation. It
follows from Lemma $3.5.b$ of \cite{Lins1980} applied to
$c_{P^\sim}$ shows that $b_P(x)=b_P(y)$ if, and only if, the
cycles $c_P(x)$ and $c_P(y)$ are homologous mod $2$.

The following proposition shows that the type of $c_P(x)$ depends
only on the parity of $i_{\overline{P}}(x)$. It is an important
tool on lacet theory.

\begin{proposition}[Parity Theorem - \cite{Lins1980}] If $M$ is a map with a single $z$-gon, then $c_P(x)$ is an
$s$-cycle in $M^t$ if and only if $|i_{\overline{P}}(x)|$ is even.
\label{prop:parityTheorem}
\end{proposition}

\comecaprova To prove the result, we define a function $c'_P$ on
the vertices of $C_M$ whose image is $CS(C_M)$. For a vertex $X$
of $C_M$, the cycle $c'_P(X)$ is defined by the edges of $C_M$
occurring once in the reentrant path which starts at $z_M(X)$ and
proceeds by using $a_{M^-},\ v_{M^-}$, and $f_M$-edges (in this
order) until it reaches another vertex of $X$, which denotes the
square to which $x$ belongs. This vertex can be $v_M(X)$, in which
case we close the path by using the $f_M$-edge linking $v_M(X)$ to
$z_M(X)$; it also can be $f_M(X)$, in which case we close the path
by using the $v_M$-edge which links $f_M(X)$ to $z_M(X)$. The fact
that $x$ is a loop in $G_P$ and a parity argument show that the
first vertex of $x$ reached by the path is not $X$, and is not
$z_M(X)$. This completes the definition of $c'_P$. Observe that
$\Psi^M_c(c'_P(X'))=c_P(x)$ for any vertex $X'$ of the square of
$M$ containing $X$, corresponding to edge $x$ of $G_M$. Since
$C_M$ is cubic, the cycle $c'_P(X)$ induces a subgraph of $C_M$
which consists of a certain number of disjoint polygons, whose set
is denoted by $\Omega$. We count the vertices of these polygons.
In the square corresponding to $x$ there are two vertices which
are vertices of a polygon in $\Omega$. If a square is met twice by
the reentrant path which defines $c'_P(X)$, then all its four
vertices are vertices of polygons in $\Omega$. If a square is met
once by the reentrant path, then three of its vertices are
vertices of a polygon in $\Omega$. Evidently, if a square is not
met by the reentrant path, none of its vertices is in a polygon in
$\Omega$. Hence, the parity of the number of vertices of the
polygons in $\Omega$, which is the same as the parity of
$|c'_P(X)|$, is equal to the parity of $|i_{\overline{P}}(x)|$.
This establishes the theorem. \terminaprova

\section{Basic results for the proof of the Main Lemma}
Let $A \in {\mathcal E}$ and $e \in E$. Define $s_e(A) = \{ f \in
A \ | \ \kappa(f)=\kappa(e)\}$.

\begin{proposition} Let $M$ be any map with a single zigzag with interlace function
$i_{\overline{P}}$. Then, for $x \in E$, $b_P(x) =
i_{\overline{P}}^2(x)+s_x(i_{\overline{P}}(x))$.
\label{prop:b-em-funcao-de-i-s}\end{proposition} \comecaprova
 $b_P(x) = c_{P^\sim} \circ c_P (x) = [c_P^2 + c_P](x) = [(i_{\overline{P}} +
\kappa_P)^2 + (i_{\overline{P}} + \kappa_P)](x)$. By expanding we
get $[i_{\overline{P}}^2 + i_{\overline{P}} \kappa_P +
\kappa_Pi_{\overline{P}} + \kappa^2 + i_{\overline{P}} +
\kappa](x) = [i_{\overline{P}}^2 + i_{\overline{P}} \kappa_P +
\kappa_Pi_{\overline{P}} + i_{\overline{P}}](x)$, since
$\kappa^2=\kappa$. If $x$ is black, then
$i_{\overline{P}}\kappa_P(x)=i_{\overline{P}}(x)$ and
$[i_{\overline{P}}^2 + i_{\overline{P}} \kappa_P +
\kappa_Pi_{\overline{P}} +
i_{\overline{P}}](x)=[i_{\overline{P}}^2 +
\kappa_Pi_{\overline{P}}
](x)=i_{\overline{P}}^2(x)+s_x(i_{\overline{P}}(x)$. If $x$ is
white, then $[i_{\overline{P}}^2 + i_{\overline{P}} \kappa_P +
\kappa_Pi_{\overline{P}} +
i_{\overline{P}}](x)=[i_{\overline{P}}^2 +
\kappa_Pi_{\overline{P}} +
i_{\overline{P}}](x)=i_{\overline{P}}^2(x)+s_x(i_{\overline{P}}(x)$.
\terminaprova

\begin{proposition} ($a$) For $x, y \in E, y \in i_{\overline{P}}^2(x)
\Leftrightarrow |i_{\overline{P}}(x) \cap i_{\overline{P}}(y)|$ is
odd. Therefore $x \in i_{\overline{P}}^2(y) \Leftrightarrow y \in
i_{\overline{P}}^2(x)$. ($b$) For $x \in E, x \in b_P(x)
\Leftrightarrow x \in i_{\overline{P}}^2(x) \Leftrightarrow x \in
{\mathcal {\mathcal O}}$. \label{prop:yInI2xIffIxIyodd}
\end{proposition} {\bf Proof:} Part $(a)$ is straightforward. For
part $(b)$ we have $b_P(x)=i^2 (x)+s_x(i_{\overline{P}}(x))$.
Since $x \notin s_x(i_{\overline{P}}(x))$ it follows that $x \in
b_P(x) \Leftrightarrow x \in i_{\overline{P}}^2(x).$ By part
$(a)$, $x \in i_{\overline{P}}^2(x) \Leftrightarrow
|i_{\overline{P}}(x) \cap
i_{\overline{P}}(x)|=|i_{\overline{P}}(x)|$ is odd
$\Leftrightarrow  x \in {\mathcal {\mathcal O}}$.

\begin{proposition} $x \in {\mathcal E}, y \in {\mathcal O} \Rightarrow b_P(x) \neq b_P(y)$.
\label{prop:oddDifferEven} \end{proposition} \comecaprova Assume
that $b_P(x)=b_P(y)$. This implies that $c_{P^\sim}
(c_P(x)+c_P(y)) = \emptyset$. Therefore $c_P(x)+c_P(y) \in
{\mathcal F}$, by Lemma $3.5.b$ of [Lins 1980] applied to
$c_{P^\sim}$. It follows that $c_P(x)+c_P(y)$ is an $s$-cycle.
This is a contradiction since, by Theorem $2$, $c_P(x)$ is an
$s$-cycle, $c_P(y)$ is an $r$-cycle and their sum must be an
$r$-cycle. \terminaprova

\begin{proposition} For $x, y \in E, y \in b_P(x) \Leftrightarrow x \in b_P(y)$.
\label{prop:bSymmetric} \end{proposition} \comecaprova Consider
the equivalence $x \in b_P(y) \Leftrightarrow [(x \in
i_{\overline{P}}^2(y)$ and $x \notin s_y(i_{\overline{P}}(y)))$ or
$(x \notin i_{\overline{P}}^2(y)$ and $x \in
s_y(i_{\overline{P}}(y)))].$ From Proposition
\ref{prop:yInI2xIffIxIyodd}, $x \in i_{\overline{P}}^2(y)
\Leftrightarrow y \in i_{\overline{P}}^2(x)$. As $x \in
i_{\overline{P}}(y) \Leftrightarrow y \in i_{\overline{P}}(x)$, $x
\in s_y(i_{\overline{P}}(y)) \Leftrightarrow y \in
s_x(i_{\overline{P}}(x))$. \terminaprova

\begin{proposition} If $P'$ is kleinian and $x, y \in {\mathcal E}, b_P(x) \ne
\emptyset \ne b_P(y)$, then $b_P(x)=b_P(y)$.
\label{prop:EvenDifZeroIguais}
\end{proposition} \comecaprova If $c_P(z)$ is an $s$-cycle for all
$z \in E$, then $C_M$ would be bipartite and $P'$ would not be
kleinian. Thus, there exists $z \in {\mathcal O} \neq \emptyset$.
Consider the non-null vectors of Im($b_P$): $b_P(x), b_P(y)$ an
$b_P(z)$. As $dim(Im(b_P))=2$ and by Proposition
\ref{prop:oddDifferEven}, either $b_P(x)=b_P(y)$ or
$b_P(z)=b_P(x)+b_P(y)$. We show that the second possibility leads
to a contradiction. This possibility implies that $b_P(x)$,
$b_P(y)$ and $b_P(z)$ are the $3$ distinct non-empty elements of
$Im(b_P)$.

As $z \in b_P(z)$ we may adjust notation and suppose that $z \in
b_P(x) \backslash b_P(y)$. Define $X = \{ x' \in E \ | \
b_P(x')=b_P(x) \}$, $Y = \{ y' \in E \ | \ b_P(y')=b_P(y) \}$ and
$Z = \{ z' \in E \ | \ b_P(z')=b_P(z) \}$. The three sets $X, Y$
and $Z$ are non-empty, since $x \in X$, $y \in Y$ and $z \in Z$.
Moreover, $X \cup Y \subseteq {\mathcal E}$, because of
Proposition \ref{prop:oddDifferEven} and $x,y \in {\mathcal E}.$
Also, $Z \subseteq {\mathcal O}$, because $z \in {\mathcal O}$.

We claim that $b_P(y) \cap Z = \emptyset$. If $z' \in b_P(y) \cap
Z$, then $z' \notin b_P(x)$. Equivalently by Proposition
\ref{prop:bSymmetric}, $x \notin b_P(z')=b_P(z)$, a contradiction
because $x \in b_P(z)$. We also claim that $b_P(y) \cap Y =
\emptyset$. If $y' \in b_P(y) \cap Y$, then $y \in b_P(y') =
b_P(y)$, contradicting Proposition \ref{prop:yInI2xIffIxIyodd}.
From the $2$ claims we may conclude that $b_P(y) \subseteq X$. Let
$x'\in b_P(y) \subseteq X$ and $x''$ be an arbitrary element of
$X$. We have that $y \in b_P(x') = b_P(x'')$ and so by Proposition
\ref{prop:bSymmetric}, $x'' \in b_P(y)$. It follows that $X
\subseteq b_P(y)$, so that $b_P(y)=X$.

Next we show that $b_P(x) = Y \cup Z$. We already know that $y, z
\in b_P(x)$. Let $y'\in Y$ and $z'\in Z$. We have $y \in b_P(x)
\Leftrightarrow x \in b_P(y)=b_P(y') \Leftrightarrow y' \in
b_P(x)$ and $z \in b_P(x) \Leftrightarrow x \in b_P(z)=b_P(z')
\Leftrightarrow z' \in b_P(x)$. Therefore, $b_P(x) \supseteq Y
\cup Z$. Let $W = \{w \in E \ | \ b(w) = \emptyset\}$, so that $E
= W \dot{\cup} X \dot{\cup} Y \dot{\cup} Z$.  If $x'\in X \cap
b_P(x)$, we get a contradiction with Proposition
\ref{prop:yInI2xIffIxIyodd}, because $x \in b_P(x')=b_P(x)$. It
follows that $b_P(x) \cap (X \cup W) = \emptyset$. Thus,
$b_P(x)\subseteq Y \cup Z$ and $b_P(x) = Y \cup Z$.

We show that $b_P(z) = X \cup Z$ as follows. We already know that
$x, z \in b_P(z)$. Let $x'\in X$ and $z'\in Z$. We have $x \in
b_P(z) \Leftrightarrow z \in b_P(x)=b_P(x') \Leftrightarrow x' \in
b_P(z)$ and $z \in b_P(z)=b_P(z') \Leftrightarrow z' \in b_P(z)$.
Therefore, $b_P(z) \supseteq X \cup Z$. As $b_P(z) \cap (Y \cup W)
= \emptyset$, we get $b_P(z) = X \cup Z = X + Z$. Also, $b_P(z) =
b_P(y) + b_P(x) = X + (Y \cup Z) = X + Y + Z$.

These two expressions for $b_P(z)$ imply that $X+Z=X+Y+Z$, or $Y =
\emptyset$, which is a contradiction. Thus, the only possibility
is $b_P(x)=b_P(y)$. \terminaprova

\begin{proposition} If $P'$ is kleinian, $x \in {\mathcal E}$ and $b_P(x) \ne \emptyset$,
then $b_P(x)={\mathcal {\mathcal O}}$.
\label{prop:ImageEvenDifZeroIsTodoOdd}
\end{proposition}
\comecaprova Let $z \in b_P(x)$. Suppose $z \in {\mathcal E}$.
Since $x \in b_P(z)$,  $b_P(x)$ and $b_P(z)$ are non-empty. By
Proposition \ref{prop:EvenDifZeroIguais}, $b_P(x)=b_P(z)$. Thus $x
\in b_P(z)=b_P(x)$, contradicting Proposition
\ref{prop:yInI2xIffIxIyodd}. Therefore $b_P(x)\subseteq {\mathcal
O}$. Let $z \in b_P(x) \subseteq {\mathcal O}$ and $z' \in
{\mathcal O} \backslash b_P(x)$. Note that $b_P(z) \neq b_P(z')$,
since $x \in b_P(z) \backslash b_P(z')$. By Proposition
\ref{prop:oddDifferEven}, $b_P(x)\neq b_P(z)$ and $b_P(x) \neq
b_P(z')$. It follows that $b_P(x), b_P(z)$ and $b_P(z')$ are the
$3$ non-null vectors in the image of $b_P$. Thus, they satisfy
$b_P(x)=b_P(z)+b_P(z')$. From this equality, since $x \in b_P(z)$
and $x \notin b_P(z')$, it follows that $x \in b_P(x)$, a
contradiction to Proposition \ref{prop:yInI2xIffIxIyodd}, because
$x \in {\mathcal E}$. So, ${\mathcal O} \backslash b_P(x)$ must be
empty, that is, $b_P(x)={\mathcal O}$. \terminaprova

{\bf Proof of the Main Lemma (Proposition {\ref{prop:teorema1}})}:
Let $z_0 \in {\mathcal O} \neq \emptyset$. Suppose that
$b_P(z_0)={\mathcal O}_0 \cup {\mathcal E}_1$, with ${\mathcal
O}_0 \subseteq {\mathcal O}$ and ${\mathcal E}_1 \subseteq
{\mathcal E}$. Note that $z_0 \in {\mathcal O}_0$. Let ${\mathcal
O}_1={\mathcal O}\backslash {\mathcal O}_0$. If both ${\mathcal
O}_1$ and ${\mathcal E}_1$ are empty, then Im($b_P$) would have
dimension $1$, and we know it is $2$. Assume that ${\mathcal E}_1
\neq \emptyset$. By the Proposition
\ref{prop:ImageEvenDifZeroIsTodoOdd}, for $x\in {\mathcal E}_1$,
then $b_P(x)={\mathcal O}$ and for $x\in {\mathcal E}_0={\mathcal
E}\backslash {\mathcal E}_1$, $b_P(x)=\emptyset$. The fourth
vector in the image of $b_P$ is $b_P(z_0+x_1) = {\mathcal O}_1
\cup {\mathcal E}_1$, for $x_1 \in {\mathcal E}_1$. Clearly,
$b_P(z_1) = {\mathcal O}_1 \cup {\mathcal E}_1, \forall z_1 \in
{\mathcal O}_1$.

If ${\mathcal E}_1=\emptyset$, then $\exists z_1 \in {\mathcal
O}_1$. Note that $\{b_P(z_0), b_P(z_1)\}$ is a basis of Im($b_P$).
Let $Z_h = \{ z \in {\mathcal O} \ | \ b_P(z)=b_P(z_h)\},\ h=1,2.$
The pair $(Z_0,Z_1)$ is a partition for ${\mathcal O}$. We claim
that $ Z_h \subseteq b_P(z_h),\ h=0,1$: indeed, as $Z_h \subseteq
{\mathcal O}$, $z \in Z_h \Rightarrow z \in b_P(z) = b_P(z_h)$.
Since $b_P(z_0) \neq b_P(z_1)$, $\exists z'_0 \in b_P(z_0)
\backslash b_P(z_1)$ or $\exists z'_1 \in b_P(z_1) \backslash
b_P(z_0)$. Note that $z'_0 \in Z_0$, if it exists, because $z'_0
\in Z_1$ leads to a contradiction: $z_1 \in b_P(z_1)$ and $z_1
\notin b_P(z'_0)$. Analogously, $z'_1 \in Z_1$, if it exists. By
replacing $z_0$ by $z'_0$ in the first case, or replacing $z_1$ by
$z'_1$ in the second and restoring the original notation we may
assume that $z_1 \notin b_P(z_0)$. For arbitrary $z'_0 \in Z_0$
and $z'_1 \in Z_1$ we claim that $z'_0 \notin b_P(z'_1)$. If $z'_0
\in b_P(z'_1)$, then $z'_0 \in b_P(z'_1) = b_P(z_1)$ and $z_1 \in
b_P(z'_0)=b_P(z_0)$ contradicting what we have established. In
consequence, $b_P(z_0) \cap Z_1$ and $b_P(z_1) \cap Z_0$ are both
empty. We can conclude that $b_P(z_h) = Z_h = {\mathcal O}_h,\
h=0,1$. \terminaprova

Let a pair of variables $\gamma_k, \delta_k \in \Z_2$ be
associated with each $k \in E$. The value of $\gamma_k$ is $1$ if
$k$ is black and is $0$ if it is white. The value of $\delta_k$ is
$h \in \{0,1\}$ if $k \in {\mathcal O}_h \cup {\mathcal E}_h$.

\begin{proposition}[Linear vectorial equations up to the Klein
bottle] Let a Gauss co\-de $\overline{P}$ be given and $E =
\{1,2,\ldots,n\}$. Each $2$-face colorable realization of
$\overline{P}$ in $\S^2$, in $\R\P^2$ or in $\S^1
\times^{\raisebox{0.5mm}{\hspace{-3.2mm}$\sim$}} \S^1$ corresponds
to a solution $(\delta_1,\ldots, \delta_n, \gamma_1, \ldots,
\gamma_n)$ of the following system of $n$ linear vectorial
equations over $\Z_2$:

\begin{center}
$k \in {\mathcal O} \Rightarrow i_{\overline{P}}^2(k) + \sum_{\ell
\in i_{\overline{P}}(k)} (1+\gamma_k+\gamma_\ell) \ell =
\sum_{\ell \in E} \delta_\ell\ell + \sum_{\ell\in {\mathcal
O}}(1+\delta_k+\delta_\ell)\ell $

$k \in {\mathcal E} \Rightarrow i_{\overline{P}}^2(k) + \sum_{\ell
\in i_{\overline{P}}(k)} (1+\gamma_k+\gamma_\ell) \ell =
\delta_k\sum_{\ell \in {\mathcal O}}\ell. $
\end{center} \label{prop:corolario1}

\end{proposition}

\comecaprova Observe that the left hand side of both equations of
Proposition \ref{prop:corolario1} yields the value of $b_P(k)$,
according to Proposition \ref{prop:b-em-funcao-de-i-s}. The value
of $\gamma_k$ is $1$ if $k$ is black and is $0$ if it is white.
The value of $\delta_k$ is $h \in \{0,1\}$ if $k \in {\mathcal
O}_h \cup {\mathcal E}_h$. Given partitions $({\mathcal
O}_0,{\mathcal O}_1)$ and $({\mathcal E}_0,{\mathcal E}_1)$ and
this interpretation of the variables $\gamma_k$, $\delta_k$ in
$\Z_2$, Proposition \ref{prop:corolario1} says the same as
Proposition \ref{prop:teorema1}. \terminaprova

{\bf Proof of Proposition \ref{prop:MainTheorem} (Main Theorem):}
 Let $\xi = (\gamma^T,\delta^T)^T$, where $\delta$ and $\gamma$ are
$n$-column vectors in $\Z_2^n$. The rows of $L$ are indexed by the
pairs $k\ell=(k,\ell) \in E^2$. The first $n$ columns of $L$
correspond to the variables $(\gamma_1, \ldots, \gamma_n)$. The
last $n$ columns of $L$ correspond to the variables $(\delta_1,
\ldots, \delta_n)$. A solution $(\gamma,\delta)$ for the Gauss
code in those surfaces satisfies, by expanding the vectorial
equations given in Proposition \ref{prop:corolario1} to their
components, to the following twelve classes of implications:

\begin{center}
\begin{tabular}{rlclclcl}
$I_1$.&$k \in {\mathcal O}$,&$\ell \in {\mathcal O},$&$\ell \in
i_{\overline{P}}^2(k)\backslash
i_{\overline{P}}(k)$&$\Rightarrow$&$0\gamma_k+0\gamma_\ell+1\delta_k+1\delta_\ell$&$=$&$0$\\
$I_2$.&$k \in {\mathcal O}$,&$\ell \in {\mathcal O},$&$\ell \in
i_{\overline{P}}^2(k)\cap
i_{\overline{P}}(k)$&$\Rightarrow$&$1\gamma_k+1\gamma_\ell+1\delta_k+1\delta_\ell$&$=$&$1$\\
$I_3$.&$k \in {\mathcal O}$,&$\ell \in {\mathcal O},$&$\ell \in
i_{\overline{P}}(k)\backslash
i_{\overline{P}}^2(k)$&$\Rightarrow$&$1\gamma_k+1\gamma_\ell+1\delta_k+1\delta_\ell$&$=$&$0$\\
$I_4$.&$k \in {\mathcal O}$,&$\ell \in {\mathcal E},$&$\ell \in
i_{\overline{P}}^2(k)\backslash
i_{\overline{P}}(k)$&$\Rightarrow$&$0\gamma_k+0\gamma_\ell+0\delta_k+1\delta_\ell$&$=$&$1$\\
$I_5$.&$k \in {\mathcal O}$,&$\ell \in {\mathcal E},$&$\ell \in
i_{\overline{P}}^2(k)\cap
i_{\overline{P}}(k)$&$\Rightarrow$&$1\gamma_k+1\gamma_\ell+0\delta_k+1\delta_\ell$&$=$&$0$\\
$I_6$.&$k \in {\mathcal O}$,&$\ell \in {\mathcal E},$&$\ell \in
i_{\overline{P}}(k)\backslash
i_{\overline{P}}^2(k)$&$\Rightarrow$&$1\gamma_k+1\gamma_\ell+0\delta_k+1\delta_\ell$&$=$&$1$\\
$I_7$.&$k \in {\mathcal E}$,&$\ell \in {\mathcal O},$&$\ell \in
i_{\overline{P}}^2(k)\backslash
i_{\overline{P}}(k)$&$\Rightarrow$&$0\gamma_k+0\gamma_\ell+1\delta_k+0\delta_\ell$&$=$&$1$\\
$I_8$.&$k \in {\mathcal E}$,&$\ell \in {\mathcal O},$&$\ell \in
i_{\overline{P}}^2(k)\cap
i_{\overline{P}}(k)$&$\Rightarrow$&$1\gamma_k+1\gamma_\ell+1\delta_k+0\delta_\ell$&$=$&$0$\\
$I_9$.&$k \in {\mathcal E}$,&$\ell \in {\mathcal O},$&$\ell \in
i_{\overline{P}}(k)\backslash
i_{\overline{P}}^2(k)$&$\Rightarrow$&$1\gamma_k+1\gamma_\ell+1\delta_k+0\delta_\ell$&$=$&$1$\\
$I_{10}$.&$k \in {\mathcal E}$,&$\ell \in {\mathcal E},$&$\ell \in
i_{\overline{P}}^2(k)\backslash
i_{\overline{P}}(k)$&$\Rightarrow$&$0\gamma_k+0\gamma_\ell+0\delta_k+0\delta_\ell$&$=$&$1$\\
$I_{11}$.&$k \in {\mathcal E}$,&$\ell \in {\mathcal E},$&$\ell \in
i_{\overline{P}}^2(k)\cap
i_{\overline{P}}(k)$&$\Rightarrow$&$1\gamma_k+1\gamma_\ell+0\delta_k+0\delta_\ell$&$=$&$0$\\
$I_{12}$.&$k \in {\mathcal E}$,&$\ell \in {\mathcal E},$&$\ell \in
i_{\overline{P}}(k)\backslash
i_{\overline{P}}^2(k)$&$\Rightarrow$&$1\gamma_k+1\gamma_\ell+0\delta_k+0\delta_\ell$&$=$&$1.$\\
\end{tabular}
\end{center}

Note that the pairs $(k,\ell)$ which do not appear in the left
hand side of the above implications are precisely the ones in
which  $\ell \notin  \left[i_{\overline{P}}(k) \cup
i_{\overline{P}}^2(k)\right]$. These pairs imply no restriction
and so it is safe to write the thirteenth class of implications
$$I_{13}. \ \ell \notin \left[i_{\overline{P}}(k) \cup
i_{\overline{P}}^2(k)\right] \Rightarrow
0\gamma_k+0\gamma_\ell+0\delta_k+0\delta_\ell=0,$$ producing a
partition of $E^2$ into 13 classes. It is now an easy matter to
display the unique matrix $L \in \Z_2^{n^2 \times 2n}$ and the
unique column vector $r \in \Z_2^{2n}$, such that the $n^2$
implications above (classified in 13 types) are equivalent to the
system $L\xi = r$. The only possibility for $(L,r)$ is to have
them implying that the $k\ell$-th equation of the system $L\xi =
r$ is $$[(L)_{k\ell,k}] \cdot \gamma_{k}+[(L)_{k\ell,\ell}] \cdot
\gamma_{\ell} + [(L)_{k\ell,k+n}]\cdot
\delta_{k}+[(L)_{k\ell,\ell+n}]\cdot \delta_{\ell} =
(r)_{k\ell}.$$ Thus, $(L)_{k\ell,p} = 0,$ if $p \notin
\{k,\ell,k+n,\ell+n\}$. Moreover, if $k\ell$ induces an
implication of type $I_q, 1 \le q \le 13$, then $(L)_{k\ell,k}$,
$(L)_{k\ell,\ell}$, $(L)_{k\ell,k+n}$, $(L)_{k\ell,\ell+n}$ and
$(r)_{k\ell}$ coincide with the coefficients in the right hand
side of implication $I_q$.

Each solution of the system $L\xi = r$  satisfies all the $n^2$
implications. Conversely, given $(\gamma,\delta)$ satisfying all
the $n^2$ implications, $\xi = (\gamma^T,\delta^T)^T$ is a
solution of $L\xi = r$. Therefore, if there is no solution for the
set of implications, the system $L\xi = r$ is inconsistent. In
this case, row operations produce a $\nu \in \Z_2^{n^2}$ such that
$\nu^T L = 0$ and $\nu^T r = 1$, giving a short proof of the
inconsistence: $L\xi = r \Rightarrow \nu^TL\xi = \nu^T r
\Rightarrow 0 = 0\xi = (\nu^T L) \xi = \nu^T(L \xi) = \nu^T r =
1.$

In fact we have used with all the $n^2$ pairs of elements of $E^2$
just for conciseness of the argument. It is easy to show that the
pairs of equal elements $(k,k)$ induce trivial restrictions and
that $(k,\ell)$ and $(\ell,k)$ induce the same restriction. Thus
we have at most $(n^2-n)/2$ restrictions. Of these, each one
coming from $k\ell$ inducing an implication in the class $I_{13}$
clearly does not need to be considered. So we need to use only $m$
of the potential $n^2$ equations, where $m \le n(n-1)/2$.
\terminaprova

As an example, the system of $n^2=8^2=64$ equations for the
example in Fig. $1$ simplifies, when the trivial and duplicated
equations are discarded, to the following system, consisting of
only 16 equations in the 16 variables $\gamma_1,\ldots,
\gamma_8,\delta_1,\ldots,\delta_8$:

\begin{center}
$ \begin{array}{|c|c|c|c|}\hline
\delta_1 = 1 + \gamma_1 + \gamma_2 & \delta_1 = 1 & \delta_1 = 1 + \gamma_1 + \gamma_6 &
\delta_1 = \gamma_1 + \gamma_7 \\ \hline
\delta_1 = \gamma_1 + \gamma_8 &  \delta_2 + \delta_6 = 0&
\delta_2 + \delta_7 = 0 &  \delta_2 + \delta_8 = 0\\ \hline
\delta_3 = \gamma_3 + \gamma_7 & \delta_3 = \gamma_3 + \gamma_8 & \delta_4 + \delta_5 = 0 &
\delta_4 + \delta_6 = \gamma_4 + \gamma_6\\ \hline
\delta_5 + \delta_6 = \gamma_5 + \gamma_6 &  \delta_6 + \delta_7 = 0 &
\delta_6 + \delta_8 = 0 & \delta_7 + \delta_8 = \gamma_7 + \gamma_8\\ \hline
\end{array}$
\end{center}

It is interesting to observe that the solution is by no means
unique. In this case, the dimension of the solution space is $4$.
The specific solution corresponding to Fig. $1$ is the second on
the left table below. All the sixteen solutions are displayed
below. As expected, the solution set is closed if we interchange
black and white vertices: $\gamma'_k=1-\gamma_k$, for all $k \in
E$.

\begin{center}
{\tiny
$\begin{array} {|c|c|c|c|c|c|c|c|c|c|c|c|c|c|c|c|} \hline
\gamma_1&\gamma_2&\gamma_3&\gamma_4&\gamma_5&\gamma_6&\gamma_7&\gamma_8&
\delta_1&\delta_2&\delta_3&\delta_4&\delta_5&\delta_6&\delta_7&\delta_8
\\ \hline \hline
0&0&1&0&0&0&1&1&1&0&0&0&0&0&0&0\\ \hline
0&0&1&1&1&0&1&1&1&0&0&1&1&0&0&0\\ \hline
0&0&0&0&0&0&1&1&1&0&1&0&0&0&0&0\\ \hline
0&0&0&1&1&0&1&1&1&0&1&1&1&0&0&0\\ \hline
0&0&1&1&1&0&1&1&1&1&0&0&0&1&1&1\\ \hline
0&0&1&0&0&0&1&1&1&1&0&1&1&1&1&1\\ \hline
0&0&0&1&1&0&1&1&1&1&1&0&0&1&1&1\\ \hline
0&0&0&0&0&0&1&1&1&1&1&1&1&1&1&1\\ \hline
\end{array}$ \hspace{5mm}
$\begin{array} {|c|c|c|c|c|c|c|c|c|c|c|c|c|c|c|c|} \hline
\gamma_1&\gamma_2&\gamma_3&\gamma_4&\gamma_5&\gamma_6&\gamma_7&\gamma_8&
\delta_1&\delta_2&\delta_3&\delta_4&\delta_5&\delta_6&\delta_7&\delta_8
\\ \hline \hline
1&1&0&1&1&1&0&0&1&0&0&0&0&0&0&0\\ \hline
1&1&0&0&0&1&0&0&1&0&0&1&1&0&0&0\\ \hline
1&1&1&1&1&1&0&0&1&0&1&0&0&0&0&0\\ \hline
1&1&1&0&0&1&0&0&1&0&1&1&1&0&0&0\\ \hline
1&1&0&0&0&1&0&0&1&1&0&0&0&1&1&1\\ \hline
1&1&0&1&1&1&0&0&1&1&0&1&1&1&1&1\\ \hline
1&1&1&0&0&1&0&0&1&1&1&0&0&1&1&1\\ \hline
1&1&1&1&1&1&0&0&1&1&1&1&1&1&1&1\\ \hline
\end{array}$ }

\end{center}

\section{2-Face colorable lacets on general surfaces}

A basic question about a Gauss code $\overline{P}$, apparently not
considered before, is to determine its {\em connectivity},
$conn(\overline{P})$, defined as the minimum connectivity among
the connectivities of the surfaces of $P'$, which realize
$\overline{P}$ as a lacet (2-face colorable or not).

Given a Gauss code $\overline{P}$, let $\gamma$ be any $0-1$
vector indexed by $E$. Vector $\gamma$ induces a map $P_\gamma$
with a single vertex as follows: the edges around the single vertex are as in
$\overline{P}$ and the orientation reversing loops are the
edges $x$ of $P_\gamma$ with $\gamma_x=1$; the others, with $\gamma_x=0$ are
orientation preserving. The phial map, $M_\gamma$, of $P_\gamma$
is a map with a single zigzag and its medial map $P'_\gamma$ is a
lacet for $\overline{P}$ on a surface $S_\gamma$ which has
$B_\gamma=\{x \in E \ | \ \gamma_x=1\}$ and $W_\gamma = \{ x\in E
\ | \ \gamma_x=0\}$ as its sets of black and white vertices. The
dimension of the image of $b_{P_\gamma}$ is the connectivity of
$S_\gamma$. The important fact is that both $b_{P_\gamma}$ and
$S_\gamma$ are defined from $\overline{P}$ and $\gamma$.

\begin{proposition}Given a Gauss code $\overline{P}$ over $E$ and an arbitrary $0-1$ vector
$\gamma$ indexed by $E$, let
$p=dim(Im(b_{P_\gamma}))=conn(S_\gamma$). Then the following
quadratic system of $n$ vectorial equations on the $2np$ variables
$\delta_{11}$, $\delta_{12}$, $\ldots$, $\delta_{1p}$, $\ldots$,
$\delta_{n1}$, $\delta_{n2}$, $\ldots$, $\delta_{np}$, $\ldots$,
$\epsilon_{11}$, $\epsilon_{12}$, $\ldots$, $\epsilon_{1p}$,
$\ldots$, $\epsilon_{n1}$, $\epsilon_{n2}$, $\ldots$,
$\epsilon_{np}$ is solvable: for each $x \in E$,
\begin{center}
$ \sum_{z \in E} (\sum_{j=1}^p \delta_{zj}\epsilon_{xj})z =
i_{\overline{P}}^2(x) + \sum_{y \in i_{\overline{P}}(x)}
(1+\gamma_x+\gamma_y) y = b_{P_\gamma}(x).$
\end{center}
\label{prop:quadraticVectorialSystem}
\end{proposition}
\comecaprova Since the right hand side of each of the above
equations is the value of $b_{P_\gamma}(x)$ it is enough to prove
that the coefficient of $z$ on the right hand side is $1$ or $0$
according to $z \in b_{P_\gamma}(x)$ or not. Henceforth, to
simplify the notation, we drop the subscripts $\gamma$. Observe
that $z \in b_P(x) = c_P(c_{P^\sim}(x))$ if and only if $| c_P(z)
\cap c_{P^\sim}(x) |$ is odd. Moreover, $c_P(z) \cap
c_{P^\sim}(x)$ and $c_P(x) \cap c_{P^\sim}(z)$ have the same
parity. The crucial observation is that the parity of $c_P(z) \cap
c_{P^\sim}(x)$ is the intersection number of the mod 2 homology
classes (Giblin$[1977]$) of the cycle $c_P(z)$ on the map $M$ and
of the cycle $c_{P^\sim}(x)$ on the dual map $D$. As we only use
homology mod 2, henceforth we drop the {\em mod 2}.

Let $N$ be an arbitrary map with a single vertex and a single face
defined on the same surface $S=S_\gamma$. $N$ is formed by loops
$v_1$, $v_2$, \ldots, $v_p$. Denote by $N^\star$ the geometrical
dual of $N$ formed by loops $v^\star_1$, $v^\star_2$, $\ldots$,
$v^\star_p$. Loop $v_j$ crosses loop $v^\star_j$ once and is
disjoint from the others. The $v_j$'s form a basis for the
homology of $S$ and the $v^\star_j$'s the dual basis. We are
considering all the maps $M$, $D$, $N$ and $N^\star$
simultaneously on the same surface $S$. Each crossing between $N$
and $D$ or between $N$ and $N^\star$ is between dual edges. All
the other crossings are irrelevant for our argument. Let $\sim$
denote the relation {\em ``is homologous to''}. Since the $v_j$'s
form a basis for the homology of $S$, there are unique scalars
$\delta_{zj}$ such that $c_P(z) \sim \sum_{j=1}^p \delta_{zj}v_j$.
Similarly, the $v^\star_j$'s form a basis for the homology of $S$.
Thus, there are unique scalars $\epsilon_{xj}$ such that
$c_{P^\sim}(x) \sim \sum_{j=1}^p \epsilon_{xj}v^\star_j$. The
crossing number between the homology classes of $c_P(z)$ and
$c_{P^\sim}(x)$ can be computed from the crossing number between
$\sum_{j=1}^p \delta_{zj}v_j$ and $\sum_{j=1}^p
\epsilon_{xj}v^\star_j$. This number is simply $\sum_{j=1}^p
\delta_{zj}\epsilon_{xj}$, proving the proposition. \terminaprova

\begin{proposition} [Vectorial quadratic equations for $P'$ on arbitrary surfaces]
Let a Gauss code $\overline{P}$ over $E$ be given. There exists a
$0-1$ vector $\gamma$ indexed by $E$ inducing a medial $P'_\gamma$
which is a lacet for $\overline{P}$ in a surface of connectivity
at most $p$ if and only if the following quadratic system of $n$
vectorial equations on the $(1+2p)n$ variables $\gamma_1$,
$\gamma_2$, \ldots, $\gamma_n$, $\delta_{11}$, $\delta_{12}$,
$\ldots$, $\delta_{1p}$, $\ldots$, $\delta_{n1}$, $\delta_{n2}$,
$\ldots$, $\delta_{np}$, $\ldots$, $\epsilon_{11}$,
$\epsilon_{12}$, $\ldots$, $\epsilon_{1p}$, $\ldots$,
$\epsilon_{n1}$, $\epsilon_{n2}$, $\ldots$, $\epsilon_{np}$ is
solvable: for each $x \in E$,
\begin{center}
$ \sum_{y \in i_{\overline{P}}(x)} (1+\gamma_x+\gamma_y) y +
\sum_{z \in E} (\sum_{j=1}^p \delta_{zj}\epsilon_{xj})z  =
i_{\overline{P}}^2(x).$
\end{center}
\label{prop:VectQuadEquations}
\end{proposition}

\comecaprova If there exists a $\gamma$ with
$dim(Im(b_{P_\gamma}))=conn(S_\gamma) = q \le p$, then apply the
previous proposition with $\overline{P}$ and this $\gamma$ as the
input parameters. We get values for $\delta_{xj}$ and
$\epsilon_{xj}$ with $x \in E$ and $j \in \{1,2,\ldots,q\}$. By
defining $\delta_{xj}=0$ and $\epsilon_{xj}=0$ for $x\in E$ and $j
\in \{q+1, \ldots, p\}$, we produce a solution of the above
system.

To prove the opposite implication, assume that the system has a
solution $(\gamma, \delta, \epsilon)$: $\gamma$ is an $n$-vector
indexed by $E$, $\delta$ and $\epsilon$ are $n \times p$ matrices.
It is enough to show that $conn(S_\gamma)=dim(Im(b_{P_\gamma}))\le
p$. Let $\{x_1,\ldots,x_p,x_{p+1}\}$ be an arbitrary subset with
$p+1$ elements of $E$ whose only restriction is that
$b_{P_\gamma}(x_i)\neq 0$, $1\le i\le p+1$. Consider the $(p+1)
\times p$ matrix $K$ whose $(i,j)$ entry is $\epsilon_{x_i,j}$.
Since $K$ has more rows than columns, there exists a subset $I
\subseteq \{1,\ldots,p,p+1\}$ such that the sum of the $K$-rows
indexed by $I$ is the zero row. Obviously $I\neq \emptyset$. Since
$(\gamma, \delta, \epsilon)$ is a solution, we get $\sum_{i \in I}
b_{P_\gamma}(x_i) = \sum_{z \in E} (\sum_{j=1}^p \delta_{zj}
(\sum_{i \in I} \epsilon_{x_i,j}))z = \sum_{z \in E} (\sum_{j=1}^p
\delta_{zj} 0)z = \sum_{z\in E} 0z = 0.$ Thus, any set of $p+1$
$b_{P_\gamma}(x_i)$'s is linearly dependent. Therefore,
$dim(Im(b_{P_\gamma}))$  is at most $p$, proving the result.
\terminaprova

Let $conn_2(\overline{P})$ denote the minimum connectivity of a
surface where $\overline{P}$ is realizable as a $2$-face colorable
lacet.

\begin{proposition} [Quadratic equations for $P'$ on arbitrary
surfaces] Given a Gauss code $\overline{P}$,
$conn_2(\overline{P})\le p$ if and only if the following system of
$n^2$ quadratic equations (on the $n(2p+1)$ variables of the
previous theorem) over $\Z_2$ is solvable: for $x,y \in E$,

\begin{center}
$\alpha_{xy}(\gamma_x+\gamma_y) + \sum_{j=1}^p
\delta_{yj}\epsilon_{xj} = \beta_{xy},$
\end{center}
where the constants $\alpha_{xy}$ and $\beta_{xy}$ are defined as
$\alpha_{xy}=1$, if $y \in i_{\overline{P}}(x)$, $\alpha_{xy}=0$,
if $y \notin i_{\overline{P}}(x)$, $\beta_{xy}=1$, if $y \in
i_{\overline{P}}(x)+i_{\overline{P}}^2(x)$ and $\beta_{xy}=0$, if
$y \notin i_{\overline{P}}(x)+i_{\overline{P}}^2(x)$.
\label{prop:QuadEquations}
\end{proposition}
\comecaprova The result follows from the previous proposition:
just expand the vectorial equations to their components.
\terminaprova

\section{Concluding Remarks}
Whether or not the above quadratic system can be efficiently
solved is an open question. Of course, if it can, then we also can
efficiently obtain $conn_2(\overline{P})$. In the small examples
we have tested, they are easy to solve via a package like MAPLE.
We have made some theoretical progress in the case of the torus.
This will be reported in a new paper currently under preparation.
Anotherr line of research is to solve the lacet problem when the
condition of face 2-colorablility is not imposed. This has been
recently achieved in the case of the projective plane
\cite{Oliveira-Lima2003}.

\section{Acknowledgements}
We thank Bruce Richter for bringing the Crapo-Rosenstiehl paper to
our attention. We also thank an anonymous referee who read
carefully a previous version of this work giving valuable
suggestions for its improvement. The first author acknowledges the
partial support of CNPq (contract number 30.1103/80).

\end{document}